# EXTENSIONS OF VECTOR BUNDLES AND RATIONALITY OF CERTAIN MODULI SPACES OF STABLE BUNDLES

WEI-PING LI AND ZHENBO QIN

ABSTRACT. In this paper, it is proved that certain stable rank-3 vector bundles can be written as extensions of line bundles and stable rank-2 bundles. As an application, we show the rationality of certain moduli spaces of stable rank-3 bundles over the projective plane $\mathbb{P}^2$.

## 1. Introduction.

Recently, there have been renewed interests in the moduli spaces of higher rank stable vector bundles over algebraic surfaces. The geometry and topology of these moduli spaces have been extensively studied. General results include the nonemptiness (see [21, 18, 13]), generic smoothness and irreducibility (see [3, 8, 27, 10, 24, 11, 25]), and variations of these moduli spaces (see [12, 14, 6, 9, 23]). There are also some specific results for these moduli spaces over particular surfaces (see [4, 5, 16, 18, 19, 26, 28] for some of the papers). In fact, quite a lot of work has been concentrated on moduli spaces of stable vector bundles over surfaces with Kodaira dimension $-\infty$. One important approach to study moduli spaces is to write the generic stable bundles as extensions of two coherent sheaves. For instance, using this approach, we proved in [18] a nonemptiness result, and showed that certain moduli spaces on geometrically ruled surfaces are unirational. However, whether these moduli spaces are rational is a delicate problem (we refer to [2, 15, 21, 7, 20] for the case of moduli spaces of stable rank-2 bundles on the projective plane $\mathbb{P}^2$).

A well known result for vector bundles over curves is that any vector bundle of rank $\geq 2$ can be written as an extension of lower rank vector bundles. For surfaces, we may not be able to get such a nice result. However, one may ask a weaker question:

Question: *Can a generic stable vector bundle of rank $\geq 3$ be written as an extension of stable bundles of lower ranks?*

If the answer to the question above is yes, we may ask further that whether the extension is unique. On the one hand, this will not be true in general since this would imply that the moduli space has a fibration over the product of two moduli spaces of lower rank stable bundles such that a generic fiber would be unirational. It follows that the Kodaira dimension of the moduli space would be $-\infty$. However, based on Jun Li's work [17] on the case of rank-2, one should expect that the

---

The first author was partially supported by the grant DAG93/94.SC22, and the second author was partially supported by a grant from the ORAU Junior Faculty Enhancement Award Program and by NSF grant DMS-94-00729.

moduli space of higher rank stable bundles over a surface of general type is also of general type. On the other hand, we shall prove in this paper that on a rational surface, a generic stable rank-3 bundle can be written as an extension of a stable rank-2 bundle and a line bundle. In fact, we expect that slight modifications of our arguments can generalize our results to arbitrary rank stable bundles on any ruled surface, and shall leave this to interested readers.

It turns out that in certain cases, the extensions for generic stable rank-3 bundles are "canonical". As an application, we can relate certain moduli spaces of stable rank-3 bundles on $\mathbb{P}^2$ to some moduli spaces of stable rank-2 bundles on $\mathbb{P}^2$. It is well-known that moduli spaces of stable rank-2 bundles on $\mathbb{P}^2$ are rational. It follows that we are able to show the rationality of some moduli spaces of stable rank-3 bundles on $\mathbb{P}^2$. More precisely, let $H$ be the divisor represented by a line in $\mathbb{P}^2$, and let $\mathfrak{M}_H(3; aH, c)$ be the moduli space of $H$-stable rank-3 bundles with first and second Chern classes $aH$ and $c$ respectively. Then we have the following.

**Theorem A.** *Let $c = (3n^2 + 9n + 4)/2$ for some positive integer $n$.*

(i) *Then a generic bundle $V \in \mathfrak{M}_H(3; 0, c)$ sits in an exact sequence*

$$0 \to \mathcal{O}_X(-nH) \to V \to W \to 0$$

*where $W$ is an $H$-stable rank-2 bundle;*

(ii) *Moreover, the moduli space $\mathfrak{M}_H(3; 0, c)$ is rational.*

**Theorem B.** *Let $c = (3n^2 + 7n + 2)/2$ or $(3n^2 + 11n + 8)/2$ for some positive integer $n$. Then the moduli space $\mathfrak{M}_H(3; H, c)$ is rational.*

We remark that from the work in [18, 26], it is known that any nonempty moduli spaces of stable bundles over $\mathbb{P}^2$ is unirational. Moreover, we have showed in [19] that the moduli space $\mathfrak{M}_H(3; H, c)$ is rational if $c$ is odd. Our Theorem B says that $\mathfrak{M}_H(3; H, c)$ is also rational for some even integers $c$. Unfortunately, Theorem B does not imply the rationality of $\mathfrak{M}_H(3; H, c)$ for an arbitrary even integer $c$.

**Acknowledgements:** The authors would like to thank the referee for valuable comments and suggestions. In particular, the referee kindly pointed out that there are other methods to prove Theorem A and Theorem B.

## 2. Stable rank-3 bundles as extensions of stable rank-2 bundles and line bundles.

One way to study vector bundles (or equivalently, locally free sheaves) over an algebraic surface $X$ is to use extensions. It works as follows. Let $V$ be a rank-$r$ vector bundle over $X$, and pick an ample divisor $H$ on $X$. Suppose that $h^2(X, V \otimes \mathcal{O}_X(nH)) = 0$ and $\chi(V \otimes \mathcal{O}_X(nH)) > 0$ for some integer $n$. For instance, this can be achieved by choosing $n$ to be sufficiently large such that $h^2(X, V \otimes \mathcal{O}_X(nH)) = 0$ by Serre's vanishing theorem and that $\chi(V \otimes \mathcal{O}_X(nH)) > 0$ by the Riemann-Roch formula. Then, $h^0(X, V \otimes \mathcal{O}_X(nH)) > 0$. It follows that we have an injection $\mathcal{O}_X(-nH) \hookrightarrow V$ which can be extended to $\mathcal{O}_X(F) \hookrightarrow V$ such that the quotient $V/\mathcal{O}_X(F)$ is a torsion free rank-$(r-1)$ sheaf. Thus, $V$ sits in an extension:

$$0 \to \mathcal{O}_X(F) \to V \to W \to 0 \tag{2.1}$$

where $W$ is a torsion free rank-$(r-1)$ sheaf. Let $W^{**}$ be the double dual of $W$. Then, $W^{**}$ is a reflexive sheaf and hence is locally free. Moreover, there is a canonical exact sequence relating $W$ and $W^{**}$:

$$0 \to W \to W^{**} \to Q \to 0. \tag{2.2}$$

If we apply the same procedure to $V^*$, we get a similar exact sequence:

$$0 \to \mathcal{O}_X(G) \to V^* \to U \to 0 \tag{2.3}$$

where $U$ is a rank-$(r-1)$ torsion free sheaf.

Suppose the rank of $V$ is bigger than two. Since $X$ has dimension two, we expect that for a generic $V$, a generic section in $H^0(X, V \otimes \mathcal{O}_X(-F))$ or $H^0(X, V^* \otimes \mathcal{O}_X(-G))$ induced by the exact sequence (2.1) or (2.3) respectively is nonvanishing everywhere. Equivalently, we expect that the sheaf $W$ in (2.1) or $U$ in (2.3) is locally free. In this section, we shall give an affirmative answer to this question when the vector bundle $V$ is a generic stable rank-3 bundle on some surface. We recall the definition of stability (in the sense of Mumford-Takemoto).

**Definition 2.4.** A rank-$r$ vector bundle $V$ on $X$ is $H$-stable if for any proper subsheaf $E$ of $V$ with rank strictly less than $r$, we have

$$\frac{c_1(E) \cdot H}{rk(E)} < \frac{c_1(V) \cdot H}{rk(V)}.$$

For $c_1 \in \text{Pic}(X)$ and $c_2 \in \mathbb{Z}$, we use $\mathfrak{M}_H(r; c_1, c_2)$ to stand for the moduli space of $H$-stable rank-$r$ vector bundles $V$ with $c_1(V) = c_1$ and $c_2(V) = c_2$.

From now on, we study $\mathfrak{M} = \mathfrak{M}_H(3; c_1, c_2)$, i.e., the moduli space of $H$-stable rank-3 bundles $V$ with $c_1(V) = c_1$ and $c_2(V) = c_2$. Our method is to stratify the moduli space $\mathfrak{M}$ into various pieces $\mathfrak{M}_{F,G}$ defined as follows.

**Definition 2.5.** For fixed divisors $F$ and $G$ on $X$, we let $\mathfrak{M}_{F,G}$ be the subset of $\mathfrak{M}$ parameterizing all the bundles $V$ such that $V$ sits in (2.1) with

$$F \cdot H = \max\{F' \cdot H | \mathcal{O}_X(F') \text{ is a subline bundle of } V\}$$

and that its dual bundle $V^*$ sits in (2.3) with

$$G \cdot H = \max\{G' \cdot H | \mathcal{O}_X(G') \text{ is a subline bundle of } V^*\}.$$

Next, we are going to estimate the dimension of each stratum $\mathfrak{M}_{F,G}$ for fixed $F$ and $G$. The key to make this work is the observation below.

**Lemma 2.6.** Let $V \in \mathfrak{M}_{F,G}$, $W$ and $U$ be in (2.1) and (2.3) respectively.
  (i) If $G \cdot H < (F - c_1) \cdot H/2$, then $W^{**}$ is $H$-stable;
  (ii) If $F \cdot H < (G + c_1) \cdot H/2$, then $U^{**}$ is $H$-stable;
  (iii) Either $W^{**}$ or $U^{**}$ is $H$-stable.

*Proof.* (i) Assume that $G \cdot H < (F - c_1) \cdot H/2$. Suppose $W^{**}$ is not $H$-stable. Let $\mathcal{O}_X(D)$ be a destablizing line bundle of $W^{**}$, i.e. we have the exact sequence

$$0 \to \mathcal{O}_X(D) \to W^{**} \to \mathcal{O}_X(c_1 - F - D) \otimes I_Z \to 0$$

with $D \cdot H \geq c_1(W^{**}) \cdot H/2 = (c_1 - F) \cdot H/2$. Then, we have a nontrivial map

$$V \to \mathcal{O}_X(c_1 - F - D) \otimes I_Z,$$

or equivalently a nontrivial map $\mathcal{O}_X(D + F - c_1) \hookrightarrow V^*$. Since $V \in \mathfrak{M}_{F,G}$, by the assumption on $G$, we have $G \cdot H \geq (D + F - c_1) \cdot H$. Thus,

$$G \cdot H \geq (D + F - c_1) \cdot H \geq \frac{1}{2}(c_1 - F) \cdot H + (F - c_1) \cdot H = \frac{1}{2}(F - c_1) \cdot H.$$

But this contradicts to our assumption that $G \cdot H < (F - c_1) \cdot H/2$.

(ii) Follows from (i) by writing $V$ as $(V^*)^*$ and switching the roles of $F$ and $G$.

(iii) In view of (i) and (ii), it suffices to show that either $G \cdot H < (F - c_1) \cdot H/2$ or $F \cdot H < (G + c_1) \cdot H/2$. Assume the contrary. Then, we would have $G \cdot H \geq (F - c_1) \cdot H/2$ and $F \cdot H \geq (G + c_1) \cdot H/2$. Thus,

$$F \cdot H \geq (G + c_1) \cdot H/2 \geq \frac{1}{4}(F - c_1) \cdot H + \frac{1}{2}c_1 \cdot H.$$

It follows that $F \cdot H \geq (c_1 \cdot H)/3$, contradicting to the $H$-stability of $V$. □

We define $\mathfrak{M}_{F^s,G}$ to be the subset of $\mathfrak{M}_{F,G}$ consisting of those $V$ whose corresponding $W^{**}$ is $H$-stable, and $\mathfrak{M}_{F,G^s}$ to be the subset of $\mathfrak{M}_{F,G}$ consisting of those $V$ whose corresponding $U^{**}$ is $H$-stable. Then, Lemma 2.6 (iii) says that $\mathfrak{M}_{F,G}$ can be split into two parts, $\mathfrak{M}_{F^s,G}$ and $\mathfrak{M}_{F,G^s}$. In the following, we shall only study $\mathfrak{M}_{F^s,G}$; by symmetry, similar results for $\mathfrak{M}_{F,G^s}$ can be obtained.

In the rest of this section, we will assume that $H \cdot K_X < 0$ where $K_X$ is the canonical divisor of $X$. Then, $X$ is necessarily a surface with Kodaira dimension $-\infty$. The following result of Maruyama (Corollary 6.7.3 in [21]) is well-known.

**Lemma 2.7.** *Assume that $H \cdot K_X < 0$. Then the moduli space $\mathfrak{M}_H(r; c_1, c_2)$ is either empty or smooth with the expected dimension*

$$2rc_2 - (r-1)c_1^2 - (r^2 - 1)\chi(\mathcal{O}_X).$$

**Lemma 2.8.** *Assume that $H \cdot K_X < 0$. If $V \in \mathfrak{M}_{F^s,G}$, then the corresponding $W$ in (2.1) varies in a family of dimension no bigger than*

$$4c_2 + 3F^2 - 2c_1 \cdot F - c_1^2 - 3\chi(\mathcal{O}_X).$$

*If equality holds, then $W$ is locally free for a generic $V \in \mathfrak{M}_{F^s,G}$.*

*Proof.* Put $c_i' = c_i(W)$ for $i = 1, 2$. Let $\overline{\mathfrak{M}_H(2; c_1', c_2')}$ be the moduli space of rank-2 Gieseker $H$-semistable sheaves with first and second Chern classes $c_1'$ and $c_2'$ respectively. From (2.1), we have $c_1' = c_1 - F$ and $c_2' = c_2 - F \cdot (c_1 - F)$. Since $V \in \mathfrak{M}_{F^s, G}$, $W^{**}$ and $W$ are $H$-stable. Thus, $W \in \overline{\mathfrak{M}_H(2; c_1', c_2')}$. By the Corollary 1.5 in [1], $\mathfrak{M}_H(2; c_1', c_2')$ is open and dense in $\overline{\mathfrak{M}_H(2; c_1', c_2')}$. Now the conclusions follow from the observation that $\mathfrak{M}_H(2; c_1', c_2')$ has dimension

$$4c_2' - (c_1')^2 - 3\chi(\mathcal{O}_X) = 4c_2 + 3F^2 - 2c_1 \cdot F - c_1^2 - 3\chi(\mathcal{O}_X). \quad \square$$

Next, we shall further assume that $X$ is a rational surface, i.e. $q(X) = 0$. Then, $H^1(X, \mathcal{O}^*)$ and $H^2(X, \mathbb{Z})$ are isomorphic. This is a technical assumption, i.e. in this case, the moduli spaces of the line bundles $\mathcal{O}_X(F)$ and $\mathcal{O}_X(G)$ in (2.1) and (2.3) are discrete sets, which makes our estimation of dimensions of various parts of moduli spaces much easier (see the proof of theorem 2.10). One may consider non-rational ruled surfaces by using the similar arguments with some care and should get a similar result as the Corollary 2.11 below.

Under this assumption, we prove the following lemma.

**Lemma 2.9.** *Assume that $H$ is an ample divisor on a rational surface $X$ with $H \cdot K_X < 0$. Let $V \in \mathfrak{M}_{F^s, G}$, and $W$ be in (2.1). Then,*
  (i) $h^0(X, W \otimes \mathcal{O}_X(K_X - F)) = 0$;
  (ii) $h^2(X, W \otimes \mathcal{O}_X(K_X - F)) = 0$;
  (iii) $h^2(X, V \otimes \mathcal{O}_X(-F)) = 0$.

*Proof.* (i) First, we claim that $H^0(X, V \otimes \mathcal{O}_X(K_X - F)) = 0$. In fact, if otherwise, we would get an injection $\mathcal{O}_X(F - K_X) \hookrightarrow V$. But $H \cdot (F - K_X) = H \cdot F - H \cdot K_X > H \cdot F$, contradicting with the assumption that $V \in \mathfrak{M}_{F^s, G}$.

Next, tensoring (2.1) by $\mathcal{O}_X(K_X - F)$ gives

$$0 \to \mathcal{O}_X(K_X) \to V \otimes \mathcal{O}_X(K_X - F) \to W \otimes \mathcal{O}_X(K_X - F) \to 0.$$

Taking cohomology of this exact sequence, we get

$$H^0(X, V \otimes \mathcal{O}_X(K_X - F)) \to H^0(X, W \otimes \mathcal{O}_X(K_X - F)) \to H^1(X, \mathcal{O}_X(K_X)).$$

Since the first and third terms are zero, so is the second term.

(ii) Suppose that $h^2(X, W \otimes \mathcal{O}_X(K_X - F)) > 0$. Then,

$$h^0(X, W^* \otimes \mathcal{O}_X(F)) = h^2(X, W^{**} \otimes \mathcal{O}_X(K_X - F)) = h^2(X, W \otimes \mathcal{O}_X(K_X - F)) > 0.$$

Thus we would have an injection $\mathcal{O}_X(-F) \hookrightarrow W^*$. Since $W^*$ is $H$-stable,

$$-F \cdot H < \frac{c_1(W^*) \cdot H}{2} = \frac{(F - c_1) \cdot H}{2},$$

or $F \cdot H > (c_1 \cdot H)/3$, contradicting to the stability of $V$.

(iii) Suppose that $h^2(X, V \otimes \mathcal{O}_X(-F)) > 0$. Then $h^0(X, V^* \otimes \mathcal{O}_X(F + K_X)) > 0$ by the Serre duality. Hence there exists an injection $\mathcal{O}_X(-F - K_X) \hookrightarrow V^*$. Since $V^*$ is $H$-stable, we must have $(-F - K_X) \cdot H < -(c_1 \cdot H)/3$. On the other hand, we have $F \cdot H < (c_1 \cdot H)/3$ by the $H$-stability of $V$ and $K_X \cdot H < 0$ by our assumption. Hence we obtain a contradiction. $\square$

Now we can prove our main result in this section.

**Theorem 2.10.** *Assume that $H$ is an ample divisor on a rational surface $X$ with $H \cdot K_X < 0$. Then, $\dim \mathfrak{M}_{F^s,G} \leq 6c_2 - 2c_1^2 - 8\chi(\mathcal{O}_X)$; moreover, if equality holds, then for a generic bundle $V \in \mathfrak{M}_{F^s,G}$, we must have*

  (i) *$W$ is locally free and $H$-stable;*
  (ii) *$h^1(X, V \otimes \mathcal{O}_X(-F)) = 0$;*
  (iii) *$h^0(X, V \otimes \mathcal{O}_X(-F))$ is positive and equal to*

$$\frac{(3F - c_1)^2 + (3F - c_1) \cdot K_X}{2} - c_2 + F \cdot (2c_1 - 3F) + 3\chi(\mathcal{O}_X).$$

*Proof.* By the Serre duality, $\text{Ext}^1(W, \mathcal{O}_X(F)) \cong H^1(X, W \otimes \mathcal{O}_X(K_X - F))$. Thus,

$$\dim \text{Ext}^1(W, \mathcal{O}_X(F)) = \frac{-(3F - c_1)^2 + (3F - c_1) \cdot K_X}{2} + c_2 + F \cdot (3F - 2c_1) - 2\chi(\mathcal{O}_X)$$

by the Riemann-Roch formula and Lemma 2.9 (i) and (ii). Similarly, by Lemma 2.9 (iii) and the Riemann-Roch formula, $h^0(X, V \otimes \mathcal{O}_X(-F))$ is equal to

$$h^1(X, V \otimes \mathcal{O}_X(-F)) + \left[ \frac{(3F - c_1)^2 + (3F - c_1) \cdot K_X}{2} - c_2 - F \cdot (3F - 2c_1) + 3\chi(\mathcal{O}_X)) \right].$$

From the exact sequence (2.1) and Lemma 2.8, we obtain

$$\begin{aligned}
\dim \mathfrak{M}_{F^s,G} \leq & \#(\text{moduli of } W) + \dim \text{Ext}^1(W, \mathcal{O}_X(F)) - h^0(X, V \otimes \mathcal{O}_X(-F)) \\
\leq & 4c_2 + 3F^2 - 2c_1 \cdot F - c_1^2 - 3\chi(\mathcal{O}_X) \\
& + \frac{-(3F - c_1)^2 + (3F - c_1) \cdot K_X}{2} + c_2 + F \cdot (3F - 2c_1) - 2\chi(\mathcal{O}_X) \\
& - \left[ \frac{(3F - c_1)^2 + (3F - c_1) \cdot K_X}{2} - c_2 - F \cdot (3F - 2c_1) + 3\chi(\mathcal{O}_X) \right] \\
\leq & 6c_2 - 2c_1^2 - 8\chi(\mathcal{O}_X).
\end{aligned}$$

If equality holds, then $\#(\text{moduli of } W)$ is equal to

$$4c_2 + F \cdot (3F - 2c_1) - c_1^2 - 3\chi(\mathcal{O}_X);$$

so by Lemma 2.8, $W$ must be a locally free sheaf for a generic $V \in \mathfrak{M}_{F^s,G}$. Moreover, $h^1(X, V \otimes \mathcal{O}_X(-F)) = 0$, and $h^0(X, V \otimes \mathcal{O}_X(-F))$ must be equal to

$$\frac{(3F - c_1)^2 + (3F - c_1) \cdot K_X}{2} - c_2 - F \cdot (3F - 2c_1) + 3\chi(\mathcal{O}_X).$$

Finally, $h^0(X, V \otimes \mathcal{O}_X(-F))$ is positive because of (2.1). $\square$

For $\mathfrak{M}_{F,G^s}$, we have similar results. Hence we obtain the corollary below.

**Corollary 2.11.** *Assume that $H$ is an ample divisor on a rational surface $X$ with $H \cdot K_X < 0$ and that the moduli space $\mathfrak{M}_H(3; c_1, c_2)$ is nonempty. Then, for a generic bundle $V \in \mathfrak{M}_H(3; c_1, c_2)$, either $V$ sits in the exact sequence (2.1) such that $W$ is locally free and $H$-stable, or $V^*$ sits in the exact sequence (2.3) such that $U$ is locally free and $H$-stable.*

## 3. Rationality of certain moduli spaces of stable rank-3 bundles on $\mathbb{P}^2$.

In this section, the surface $X$ is the projective plane $\mathbb{P}^2$, and the ample divisor $H$ is the divisor class represented by a line. For any rank-3 bundle $V$ on $X$, by taking its dual and/or tensoring it by a suitable line bundle, we can assume that $c_1(V) = aH$ with $a = 0$ or $1$. In [4], it has been showed that the moduli space $\mathfrak{M}_H(3; 0, c)$ is nonempty if and only if $c \geq 3$, and that $\mathfrak{M}_H(3; H, c)$ is nonempty if and only if $c \geq 2$. It is also known that whenever $\mathfrak{M}_H(3; aH, c)$ is nonempty, it is irreducible (see [5, 4]) and unirational (see [19, 26]). In [19], we also proved that $\mathfrak{M}_H(3; H, c)$ is rational if $c$ is odd. In this section, as an application of our results in section 2, we shall prove that generic bundles $V$ in $\mathfrak{M}_H(3; aH, c)$ sit in the two exact sequences (3.6) and (3.7). Then we determine the rationalities of $\mathfrak{M}_H(3; 0, c)$ and $\mathfrak{M}_H(3; H, c)$ for some special values of $c$.

Assume that $c \geq (a^2 + 3a + 3)$. Then $\mathfrak{M}_H(3; aH, c)$ is nonempty by [4]. Let $V \in \mathfrak{M}_H(3; aH, c)$. Choose $n$ to be the least nonnegative integer such that

$$f(n) = \frac{3n^2 + (2a+9)n + (a^2 + 3a + 6)}{2} \geq (c+1),$$

and $m$ to be the least nonnegative integer such that

$$g(m) = \frac{3m^2 + (-2a+9)m + (a^2 - 3a + 6)}{2} \geq (c+1).$$

Since $(c+1) > f(n-1) \geq g(n-1)$ and $g(n+a) \geq f(n) \geq (c+1)$, we see that $n \leq m \leq (n+a)$. Since $g(a) = (a^2 + 3a + 3) < (c+1)$, we must have $m > a$ and $n > 0$. Therefore, we obtain the following inequalities.

**Lemma 3.1.** $n > (m-a)/2$ and $m > (n+a)/2$. □

Next, applying the Riemann-Roch formula, we obtain

$$\chi(V(n)) = f(n) = \frac{3n^2 + (2a+9)n + (a^2 + 3a + 6)}{2} - c \geq 1. \tag{3.2}$$

where we have used $V(n)$ to stand for $V \otimes \mathcal{O}_X(nH)$. Since $V$ is stable,

$$h^2(X, V(n)) = h^0(X, V^*(-n-3)) = 0.$$

Thus, $h^0(X, V(n)) > 0$ and the exact sequence (2.1) becomes

$$0 \to \mathcal{O}_X((k-n)H) \to V \to W \to 0 \tag{3.3}$$

where $0 \leq k < (a/3 + n)$. Similarly, the exact sequence (2.3) becomes

$$0 \to \mathcal{O}_X((l-m)H) \to V^* \to U \to 0 \tag{3.4}$$

where $0 \leq l < (-a/3 + m)$. Thus, increasing the integers $k$ and $l$ if necessary, we have $V \in \mathfrak{M}_{(k-n)H,(l-m)H}$ for some $k$ and $l$ with $0 \leq k < (a/3 + n)$ and $0 \leq l < (-a/3 + m)$. Since $\mathfrak{M}_H(3; aH, c)$ is irreducible and $\text{Pic}(\mathbb{P}^2)$ is discret, we conclude that there exists a unique pair of integers $(k, l)$ with $0 \leq k < (a/3+n)$ and $0 \leq l < (-a/3 + m)$ such that $\mathfrak{M}_{(k-n)H,(l-m)H}$ is open and dense in $\mathfrak{M}_H(3; aH, c)$. Now we are ready to prove our next result.

**Theorem 3.5.** *Assume that $c \geq (a^2 + 3a + 3)$ where $a = 0$ or $1$. Let $n$ be the least nonnegative integer with $\left[3n^2 + (2a+9)n + (a^2+3a+6)\right]/2 \geq (c+1)$, and $m$ be the least nonnegative integer with $\left[3m^2 + (-2a+9)m + (a^2-3a+6)\right]/2 \geq (c+1)$. Then, a generic bundle $V \in \mathfrak{M}_H(3; aH, c)$ sits in the following two exact sequences*

$$0 \to \mathcal{O}_X(-nH) \to V \to W \to 0 \tag{3.6}$$

$$0 \to \mathcal{O}_X(-mH) \to V^* \to U \to 0 \tag{3.7}$$

*where $W$ and $U$ are $H$-stable rank-2 vector bundles.*

*Proof.* We have showed that there exists a unique pair of integers $(k, l)$ with $0 \leq k < (a/3+n)$ and $0 \leq l < (-a/3+m)$ such that a generic bundle $V \in \mathfrak{M}_H(3; aH, c)$ is contained in $\mathfrak{M}_{(k-n)H,(l-m)H}$. Since $\mathfrak{M}_{(k-n)H,(l-m)H}$ is the union of two subsets $\mathfrak{M}_{(k-n)H^s,(l-m)H}$ and $\mathfrak{M}_{(k-n)H,(l-m)H^s}$, there are two possibilities.

**Case 1:** A generic bundle $V \in \mathfrak{M}_{(k-n)H,(l-m)H}$ is contained in $\mathfrak{M}_{(k-n)H^s,(l-m)H}$. Then by Theorem 2.10 (iii), for a generic bundle $V \in \mathfrak{M}_{(k-n)H^s,(l-m)H}$, we have

$$\begin{aligned}1 &\leq h^0(X, V(-(k-n)H)) \\ &= \frac{3(n-k)^2 + (2a+9)(n-k) + (a^2+3a+6)}{2} - c.\end{aligned}$$

Since $(n-k) \geq 0$, we must have $(n-k) \geq n$ by the choice of $n$; thus $k \leq 0$. It follows that $k = 0$ and that a generic bundle $V$ in the moduli space $\mathfrak{M}_H(3; aH, c)$ sits in the following two exact sequences:

$$0 \to \mathcal{O}_X(-nH) \to V \to W \to 0$$

$$0 \to \mathcal{O}_X((l-m)H) \to V^* \to U \to 0.$$

Moreover, by Theorem 2.10 (i), $W$ is an $H$-stable rank-2 bundle.

Next, by Lemma 2.6 (ii) and Lemma 3.1, $U^{**}$ is stable for a generic bundle $V \in \mathfrak{M}_{-nH^s,(l-m)H}$. It follows that a generic bundle $V \in \mathfrak{M}_{-nH,(l-m)H}$ is also contained in $\mathfrak{M}_{-nH,(l-m)H^s}$. As in the preceding paragraph, we conclude that $l = 0$ and that for a generic bundle $V \in \mathfrak{M}_{-nH,-mH}$, $U$ is an $H$-stable rank-2 vector bundle.

**Case 2:** A generic bundle $V \in \mathfrak{M}_{(k-n)H,(l-m)H}$ is contained in $\mathfrak{M}_{(k-n)H,(l-m)H^s}$. As in Case 1, we conclude that $l = 0$ and that a generic bundle $V \in \mathfrak{M}_H(3; aH, c)$ sits in the following two exact sequences:

$$0 \to \mathcal{O}_X((k-n)H) \to V \to W \to 0$$

$$0 \to \mathcal{O}_X(-mH) \to V^* \to U \to 0$$

where $U$ is an $H$-stable rank-2 vector bundle.

Next, by Lemma 2.6 (i) and Lemma 3.1, $W^{**}$ is stable for a generic bundle $V \in \mathfrak{M}_{(k-n)H,-mH^s}$. It follows that a generic bundle $V \in \mathfrak{M}_{(k-n)H,-mH}$ is also contained in $\mathfrak{M}_{(k-n)H^s, -mH}$. Hence $k = 0$, and $W$ is an $H$-stable rank-2 bundle for a generic bundle $V \in \mathfrak{M}_{-nH, -mH}$. □

**Remark 3.8.** By [4], the moduli space $\mathfrak{M}_H(3; 0, c)$ is nonempty if and only if $c \geq 3$. Since $a = 0$, we have $n = m$. Thus, if the moduli space $\mathfrak{M}_H(3; 0, c)$ is nonempty, then a generic bundle $V \in \mathfrak{M}_H(3; 0, c)$ sits in an exact sequence

$$0 \to \mathcal{O}_X(-nH) \to V \to W \to 0 \qquad (3.9)$$

where $W$ is an $H$-stable rank-2 vector bundle and $n$ is the least nonnegative integer such that $(3n^2 + 9n + 6)/2 \geq (c+1)$, or equivalently, $(3n^2 + 9n)/2 \geq (c-2)$. Moreover, we have $h^0(X, V(n)) = (3n^2 + 9n + 6)/2 - c$.

Finally, as applications of Theorem 3.5, we now prove the rationality of the moduli spaces $\mathfrak{M}_H(3; 0, c)$ and $\mathfrak{M}_H(3; H, c)$ for some values of $c$.

**Corollary 3.10.** Let $c = (3n^2 + 9n + 4)/2$ for some positive integer $n$.
 (i) Then a generic bundle $V \in \mathfrak{M}_H(3; 0, c)$ sits in an exact sequence

$$0 \to \mathcal{O}_X(-nH) \to V \to W \to 0 \qquad (3.11)$$

 where $W$ is an $H$-stable rank-2 bundle;
 (ii) Moreover, the moduli space $\mathfrak{M}_H(3; 0, c)$ is rational.

*Proof.* Since $n$ is the least nonnegative integer $t$ such that

$$(3t^2 + 9t + 6)/2 \geq (c+1),$$

(i) follows immediately from Remark 3.8. Moreover, by Remark 3.8, we see that $h^0(X, V(n)) = 1$ for a generic bundle $V \in \mathfrak{M}_H(3; 0, c)$. Thus for a generic bundle $V \in \mathfrak{M}_H(3; 0, c)$, $W$ is determined uniquely by $V$ and the extension (3.11) is canonical. Next, we notice that $W \in \mathfrak{M}_H(2; nH, c + n^2)$ and that $W$ varies in an open dense subset $\mathcal{U}$ of $\mathfrak{M}_H(2; nH, c + n^2)$ as the generic bundle $V$ varies in $\mathfrak{M}_H(3; 0, c)$. Using a standard method and the canonical extension (3.11), we can construct an algebraic projective bundle $\mathcal{P}$ over $\mathcal{U}$ such that a dense open subset of $\mathcal{P}$ is isomorphic to a dense open subset of $\mathfrak{M}_H(3; 0, c)$. Since the moduli space $\mathfrak{M}_H(2; nH, c + n^2)$ is rational [2, 15, 21, 7, 20], so is $\mathfrak{M}_H(3; 0, c)$. □

**Corollary 3.12.** Let $c = (3n^2 + 7n + 2)/2$ or $(3n^2 + 11n + 8)/2$ for some positive integer $n$. Then the moduli space $\mathfrak{M}_H(3; H, c)$ is rational.

*Proof.* Notice that $c \geq 7$. Thus, $\mathfrak{M}_H(3; H, c)$ is nonempty by [4]. Let $a = 1$ in Theorem 3.5. When $c = (3n^2 + 11n + 8)/2$, the integer $n$ here coincides with the integer $n$ in Theorem 3.5; thus, a generic $V \in \mathfrak{M}_H(3; H, c)$ sits in an exact sequence

$$0 \to \mathcal{O}_X(-nH) \to V \to W \to 0$$

where $W \in \mathfrak{M}_H(2; n+1, c + n^2 + n)$. Moreover, $h^0(X, V(n)) = 1$ by Theorem 2.10 (iii). As in the proof of Corollary 3.10, we see that $\mathfrak{M}_H(3; H, c)$ is rational.

Next, when $c = (3n^2 + 7n + 2)/2$, the integer $n$ here coincides with the integer $m$ in Theorem 3.5; thus, a generic bundle $V$ in $\mathfrak{M}_H(3; H, c)$ sits in an exact sequence

$$0 \to \mathcal{O}_X(-nH) \to V^* \to U \to 0$$

where $U$ is an $H$-stable rank-2 vector bundle and $\text{Hom}(\mathcal{O}_X(-nH), V^*)$ has dimension one. Again, we conclude that the moduli space $\mathfrak{M}_H(3; H, c)$ is rational. $\square$

**Remark 3.13.** We notice that $c = (3n^2 + 11n + 8)/2$ is even when $n \equiv 0 \pmod{4}$ and $c = (3n^2 + 7n + 2)/2$ is even when $n \equiv 2 \pmod{4}$. Hence Corollary 3.12 proves the rationality of $\mathfrak{M}_H(3; H, c)$ for some even integer $c$. But it does not imply the rationality of $\mathfrak{M}_H(3; H, c)$ for an arbitrary even integer $c$. For instance, 8 can not be written as $(3n^2 + 7n + 2)/2$ or $(3n^2 + 11n + 8)/2$ for any positive integer $n$.

DEPARTMENT OF MATHEMATICS, HKUST, CLEAR WATER BAY, KOWLOON, HONG KONG
*E-mail address*: mawpli@masu1.ust.hk

DEPARTMENT OF MATHEMATICS, OKLAHOMA STATE UNIVERSITY, STILLWATER, OK 74078, USA
*E-mail address*: zq@math.okstate.edu